\documentclass[10pt]{article}

\usepackage{graphics}
\usepackage{graphicx}

\usepackage[a4paper, left=35mm,right=35mm,top=34mm,bottom=34mm]{geometry}
\usepackage[utf8]{inputenc}
\usepackage[T1]{fontenc}
\usepackage[english]{babel}

\usepackage{enumerate}
\usepackage{graphicx}
\usepackage{hyperref}
\hypersetup{
    colorlinks=true,
    linkcolor=blue,
    filecolor=magenta,      
    urlcolor=cyan,
}

\usepackage{mathtools,amsthm,amssymb,amsfonts}
\usepackage{algorithm}
\usepackage{algorithmic}
\makeatother
\theoremstyle{plain}

\theoremstyle{definition}

\theoremstyle{remark}

\usepackage{caption} 
\usepackage[tight,footnotesize]{subfigure}

\usepackage{fancyhdr}
\usepackage{array}
\usepackage[dvipsnames]{xcolor}

\usepackage{ecp/ecpmashpc}
\usepackage{ecp/matrix_tikz/mrgcommon}
\usepackage{ecp/matrix_tikz/matrix_subdomain}

\author{
  {\normalsize Abal-Kassim Cheik Ahamed}\thanks{CUDA Research Center, \'Ecole Centrale Paris, France.}
	\and
  {\normalsize Fr\'ed\'eric Magoul\`es}\thanks{CUDA Research Center, \'Ecole Centrale Paris, France
    (correspondence, frederic.magoules@hotmail.com).}
		}	
\title{Parallel Sub-Structuring Methods for solving Sparse Linear Systems on a cluster of GPU}
\date{}

\begin{document}
\maketitle
\thispagestyle{fancy}

\begin{abstract}
\noindent The main objective of this work consists in analyzing sub-structuring method for the parallel solution of sparse linear systems with matrices arising from the discretization of partial differential equations such as finite element, finite volume and finite difference.
With the success encountered by the general-purpose processing on graphics processing units (GPGPU), we develop an hybrid multiGPUs and CPUs sub-structuring algorithm.
GPU computing, with CUDA, is used to accelerate the operations performed on each processor.
Numerical experiments have been performed on a set of matrices arising from engineering problems.
We compare C+MPI implementation on classical CPU cluster with C+MPI+CUDA on a cluster of GPU. The performance comparison shows a speed-up for the sub-structuring method up to 19 times in double precision by using CUDA.
\end{abstract}


\section{Introduction}
\label{sec:introduction}

Many engineering problems lead to the computation of large size sparse linear systems arising from the discretization of numerical methods such as finite element, finite volume and finite difference. Iterative Krylov methods are suitable to solve these kind of problems. According to the properties of the matrices involved in the linear systems, Krylov methods do not have the same efficiency~\cite{saad_iterative_2003}. In this paper, we consider problems that leads to symmetric and positive-definite matrices, which therefore places the Conjugate Gradient (CG) method as a suitable and efficient Krylov method. This method requires the computation of linear algebra operations such as scalar product, addition of vectors, matrix-vector multiplication (SpMV) that can be costly in terms of computations on a conventional processor (CPU). The use of GPGPU model can accelerate these operations. GPGPU is more efficient for large size problems. Unfortunately, most graphics cards have very limited memory size.
We therefore propose to parallelize the CG algorithm~\cite{PARALLEL:PLP:1997,PARALLEL:CZ:2002,PARALLEL:AHR:2010} using sub-structuring approach, which corresponds to a natural parallelization methodology and is easy to implement. Then, GPU is used to accelerate the computation of local linear algebra operations. Each sub-structure is associated with a single processor (CPU) and an accelerated device (GPU). 

This paper is organized as follows. The first section presents the partitionning of the data. The next section (\ref{sec:conjugate_gradient}) presents the main points of the implementation of the conjugate gradient in parallel for different partitionning, leading to the design of linear algebra operations. Section~\ref{sec:exp_evaluation} collects and analyzes the numerical results and analyze them. Matrices and experiments hardware are also presented. Finally, concluding remarks are given in Section~\ref{sec:conclusion}.

\section{GPU: An Effective Accelerator}
\label{sec:gpu_computing_csr}

\subsection{GPU computing}
\label{subsec:gpu_computing}

Since a few decades, Graphics Processing Units (GPUs) are used to accelerate scientific computation by graphics card hardware. These graphics cards were before used only for graphics applications such as Graphical User Interface. The GPU is a processor with wide computational resources. The rapid improvement of GPU performances has allowed to give to GPU Computing an important place in scientific computing, helped by the flexibility of programming on GPU with language such as CUDA.
GPU Computing or GPGPU has become essential in scientific computing when we deal with time consuming of a numerical simulation.  
Current GPUs are enable to compute simultaneously similar operations by more than one million of threads. The last graphics card of Kepler family, K40, reaches 4.29 teraflops single-precision and 1.43 teraflops double-precision peak floating point performance. The storage memory of the graphics cards have also seen their memory storage evolved in recent years, up to 12 GB for the K40. The performance of GPU algorithms strongly depends on both the configuration of the distribution of the threads on the grid~\cite{cheikahamed:2012:inproceedings-2} and the memory~\cite{Hong:2009:AMG}.
As indicated in~\cite{cheikahamed:2012:inproceedings-2}, the performance of sparse matrix-vector multiplication, which is a time consuming linear algebra operation, depends on the structure of the matrix, i.e., the pattern of non-zero value, and the format of the matrix storage. References~\cite{Lee:2014:BCA}~\cite{cheikahamed:2012:inproceedings-1}~\cite{cheikahamed:2013:inproceedings-3} confirm the influence of the distribution of the threads, the matrix structure and storage format, when solving linear systems with iterative Krylov methods on GPU.
 
\subsection{Sparse matrix formats}
\label{subsec:sparse-matrix_formats}

\begin{figure}[!ht]
\centering
\includegraphics[scale=0.25]{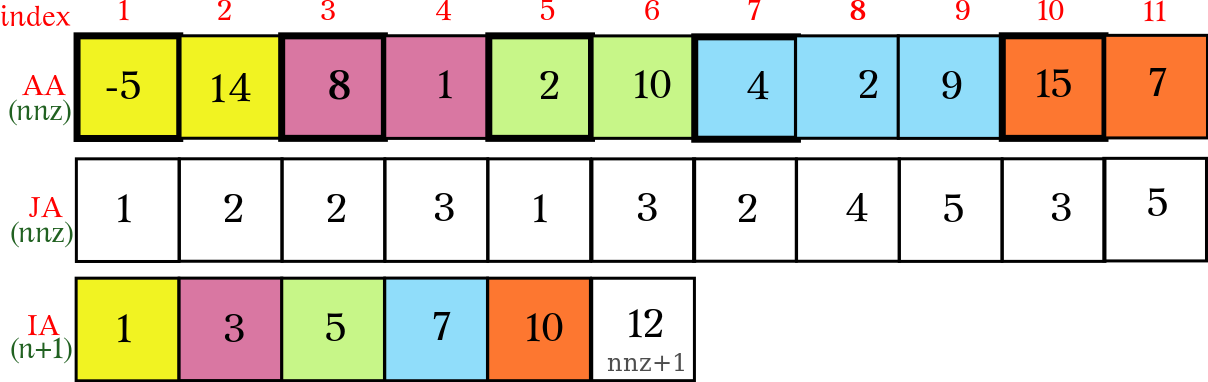}
\caption{Compressed sparse row storage format (CSR) of matrix Table~\ref{tab:matrix_and_pattern}}
\label{fig:img:csr_format}
\end{figure}

Usually, solving partial differential equations by numerical methods such as finite element method (FEM) lead to large and sparse matrix, i.e., only a few elements of the matrix are nonzero. The distribution of non-zero coefficients depends on the features of the original problem. Sparse matrix is called \emph{structured} when the non-zero values form a regular pattern along diagonals, otherwise it is called \emph{unstructured}. The performance of the algorithms strongly depends on the structure of the sparse matrix~\cite{GPU:BG:2009,GPU:BG:2008,GPU:BFGS:2003}. In terms of memory storage, sparse matrices are stored in compressed formats, which consists in only allowing memory to their non-zero coefficients. Different data storage structures exist~\cite{saad_iterative_2003} such as Compressed-Sparse Row (CSR), Coordinate (COO), ELLPACK (ELL), Hybrid (HYB), etc. In this work, we consider the CSR format. The CSR format stores the matrix using three one-dimensional arrays, as drawn in Fig.~\ref{fig:img:csr_format}. Two arrays of size $nnz$, $AA$ and $JA$ store respectively the non-zero coefficients of the matrix in consecutive rows and the column indices, i.e., $JA(k)$ is the column index in $A$ matrix of $AA(k)$. The third array, $IA$, of size $n+1$, stores pointers to the beginning of each row. $IA(i)$ and $IA(i + 1) - 1$ correspond to the beginning and the end of the $i-th$ row in arrays $AA$ and $JA$, i.e., $IA(n + 1) = nnz + 1$. An example in CSR format of matrix $A$ (Table~\ref{tab:matrix_and_pattern}) is given in Figure~\ref{fig:img:csr_format}.
\begin{table}
\centering
\begin{tabular}{m{6cm}m{4cm}}
\scalebox{1.3}{
$A=\begin{pmatrix}
\colorbox{gray}{-5} & \colorbox{gray}{14} & 0 & 0 & 0 \\
 0 & \colorbox{gray}{8} & \colorbox{gray}{1} & 0 & 0 \\
 \colorbox{gray}{2} & 0 & \colorbox{gray}{10} & 0 & 0 \\
 0 & \colorbox{gray}{4} & 0 & \colorbox{gray}{2} & \colorbox{gray}{9} \\
 0 & 0 & \colorbox{gray}{15} & 0 & \colorbox{gray}{7}
\end{pmatrix}$
}
&
\scalebox{0.7}{
\SetDisplayColumnIndexD{1}{0}
\SetDisplayRowIndexD{1}{0}
\def\nnzcoef{1/1,1/2,2/2,2/3,3/1,3/3,4/2,4/4,4/5,5/3,5/5}
\DrawGivenMatrix{(0,0)}{5}{5}{\nnzcoef}
}\\
\end{tabular}
\caption{Left (matrix), Right (matrix pattern)}
\label{tab:matrix_and_pattern}
\end{table}
%

\section{Graph and Matrix Partitionning}
\label{sec:graph_matrix_partitionning}

Undoubtedly, parallelism is the future of computing. The main step in parallel processing consists in distributing the data on the cluster processors, what is commonly called parallel distributed computing. In this section we describe how data are distributed among processors for different splitting strategies: band-row, band-column, and sub-structuring splitting.
The distribution of data is accomplished as a pre-processing step, independently to the solver code. The data such as matrix, right hand-size, vector solution and local to global, are written into file, and will be input by the solver presented in section~\ref{sec:conjugate_gradient}.

\subsection{Band-row splitting}
\label{subsec:band-row_splitting}

The partition of the set of equations leads to allocate to each processor a band of rows, corresponding to the block of the vectors it treats. In Table~\ref{tab:band-row_splitting}, these terms are located on a colorful area. The band-row splitting approach consists in partitionning the matrix $A$ of size $n \times n$ into horizontal band matrices. Each processor is in charge of the management of a band-row matrix of size $N_p \times n$ and the associated unknown vector $x$ of size $N_p \times 1$, as drawn in Table~\ref{tab:band-row_splitting}.
This method of partitionning by band-row allows to exhibit a sufficient degree of parallelism properly balanced.
For this, it implies to assign to all processors, a block of rows of same size, containing approximately the same number of non-zero coefficients. It unfortunately suffers from a major lack of granularity for implementation on a distributed memory system.
\begin{table}[!ht]
\SetDisplayIndex{1}
\def\procSy{3}
\def\procEy{6}
\def\procSa{3}
\def\procEa{6}
\def\procSx{3}
\def\procEx{6}
\centering
\setlength{\tabcolsep}{0.0pt}
\begin{tabular}{cccc}
$\;\;\;\;$ & $\;\;\;\;\;\;\;=$ & $\;$ & $\;\;\;\;$ \\
\scalebox{0.45}{
\SetSplitType{0}
\def\nnzcoef{1/1,2/1,3/1,4/1,5/1,6/1,7/1,8/1,9/1,10/1}
\SetDisplayInfoSize{1}
\SetDisplayInfoSizeProc{1}
\DrawGivenMatrixProcessor{(0,0)}{10}{1}{\nnzcoef}{\procSy}{\procEy}{\large{$n$}}{\large{$N_p$}}
}
& &
\scalebox{0.45}{
\SetSplitType{0}
\def\nnzcoef{
1/1,1/2,1/4,1/6,
2/1,2/2,2/3,2/5,2/8,
3/2,3/3,3/4,3/5,
4/1,4/3,4/4,4/7,4/8,4/10,
5/2,5/3,5/5,5/7,
6/1,6/6,6/7,6/9,
7/4,7/5,7/6,7/7,7/10,
8/2,8/4,8/8,8/10,
9/6,9/9,
10/4,10/7,10/8,10/10}
\SetDisplayInfoSize{0}
\SetDisplayInfoSizeProc{0}
\DrawGivenMatrixProcessor{(0,0)}{10}{10}{\nnzcoef}{\procSa}{\procEa}{\large{$n$}}{\large{$N_p$}}
}
&
\scalebox{0.45}{
\SetSplitType{0}
\def\nnzcoef{1/1,2/1,3/1,4/1,5/1,6/1,7/1,8/1,9/1,10/1}
\SetDisplayInfoSize{1}
\SetDisplayInfoSizeProc{1}
\DrawGivenMatrixProcessor{(0,0)}{10}{1}{\nnzcoef}{\procSx}{\procEx}{\large{$n$}}{\large{$N_p$}}
}\\
$\;\;\;Y$ & $\;\;\;\;\;\;\;$ & $A$ & $\;\;\;X$
\end{tabular}
\caption{A band-row splitting of matrix $A$}
\label{tab:band-row_splitting}
\end{table}
The algorithm of band-row partitionning of a given matrix for a given processor $p$ is presented in Algorithm~\ref{algo:band-row_partitioning}.
\begin{algorithm}[!ht]
\footnotesize
\caption{Band-row partitionning for the $p-th$ band}
\label{algo:band-row_partitioning}
\begin{algorithmic}[1]
\REQUIRE $A(n \times n, rowptr, colidx, val)$: csr matrix
\REQUIRE $ 0 < N_{b} <= n$: number of row-bands
\REQUIRE $ 1 < p <= N_{b}$: processor number ( or band number)
\STATE $numb\_rows\_per\_band \leftarrow n / N_{b}$
\STATE $N_p \leftarrow numb\_rows\_per\_band$
\IF{$p = N_{b} - 1$}
\STATE $N_p \leftarrow n - p \times numb\_rows\_per\_band$
\ENDIF
\STATE $row_{displ} \leftarrow i \times numb\_rows\_per\_band$
\STATE $nnz_{displ} \leftarrow A.rowptr[row_{displ}]$
\STATE $nnz_{count}$ $\leftarrow$ $A.rowptr[row_{displ}+N_p]$-$A.rowptr[row_{displ}]$
\STATE \COMMENT{\bf Allocate band-row matrix $A_{band}(N_p \times n, nnz_{count})$}
\STATE \COMMENT{\bf // -- copy row indices}
\STATE $A_{band}.rowptr[:] \leftarrow A.rowptr[row_{displ} + :]$
\STATE \COMMENT{\bf // -- shift row indices to local}
\STATE $A_{band}.rowptr[:] \leftarrow A_{band}.rowptr[:] - A_{band}.rowptr[1]$
\STATE $A_{band}.rowptr[N_p+1] \leftarrow nnz_{count}$
\STATE \COMMENT{\bf // -- copy cols numb and coef}
\STATE $A_{band}.rowptr[:] \leftarrow A.[nnz_{displ} + :]$
\STATE $A_{band}.val[:] \leftarrow A.[nnz_{displ} + :]$
\end{algorithmic}
\end{algorithm}

\subsection{Band-column splitting}
\label{subsec:band-column_splitting}

As for the band-row splitting, the band-column approach consists in partitionning the matrix $A$ into vertical band matrices. Each processor is in charge of the management of a band-column matrix of size $n \times N_p$. The associated unknown vector $x$ of size $N_p \times 1$ is splitting into horizontal band vectors as in band-row splitting as described in Table~\ref{tab:band-column_splitting}.
\begin{table}[!ht]
\SetDisplayIndex{1}
\def\procSy{3}
\def\procEy{6}
\def\procSa{3}
\def\procEa{6}
\def\procSx{3}
\def\procEx{6}
\centering
\setlength{\tabcolsep}{0.0pt}
\begin{tabular}{cccc}
$\;\;\;\;$ & $\;\;\;\;\;\;\;=$ & $\;$ & $\;\;\;\;$ \\
\scalebox{0.45}{
\SetSplitType{0}
\def\nnzcoef{1/1,2/1,3/1,4/1,5/1,6/1,7/1,8/1,9/1,10/1}
\SetDisplayInfoSize{1}
\SetDisplayInfoSizeProc{1}
\DrawGivenMatrixProcessor{(0,0)}{10}{1}{\nnzcoef}{\procSy}{\procEy}{\large{$n$}}{\large{$N_p$}}
}
& &
\scalebox{0.45}{
\SetSplitType{1}
\def\nnzcoef{
1/1,1/2,1/4,1/6,
2/1,2/2,2/3,2/5,2/8,
3/2,3/3,3/4,3/5,
4/1,4/3,4/4,4/7,4/8,4/10,
5/2,5/3,5/5,5/7,
6/1,6/6,6/7,6/9,
7/4,7/5,7/6,7/7,7/10,
8/2,8/4,8/8,8/10,
9/6,9/9,
10/4,10/7,10/8,10/10}
\SetDisplayInfoSize{0}
\SetDisplayInfoSizeProc{0}
\DrawGivenMatrixProcessor{(0,0)}{10}{10}{\nnzcoef}{\procSa}{\procEa}{\large{$n$}}{\large{$N_p$}}
}
&
\scalebox{0.45}{
\SetSplitType{0}
\def\nnzcoef{1/1,2/1,3/1,4/1,5/1,6/1,7/1,8/1,9/1,10/1}
\SetDisplayInfoSize{1}
\SetDisplayInfoSizeProc{1}
\DrawGivenMatrixProcessor{(0,0)}{10}{1}{\nnzcoef}{\procSx}{\procEx}{\large{$n$}}{\large{$N_p$}}
}\\
$\;\;\;Y$ & $\;\;\;\;\;\;\;$ & $A$ & $\;\;\;X$
\end{tabular}
\caption{A band-column splitting of matrix $A$}
\label{tab:band-column_splitting}
\end{table}
Algorithm~\ref{algo:band-column_partitioning} describes the column-row partitionning procedure of a given matrix for a given processor $p$. According to the structure of CSR format, the computation of the number of non-zero values of each band requires a particular calculation, unlike the row partitionning. At line 8 of the Algorithm~\ref{algo:band-column_partitioning}, we recover the number of non-zero values computed outside the routine. All non-zero values of all processors are stored into an independant array, which is built using the same test process described at line 14 of Algorithm~\ref{algo:band-column_partitioning}.
\begin{algorithm}[!ht]
\footnotesize
\caption{Band-column partitionning for the $p-th$ band}
\label{algo:band-column_partitioning}
\begin{algorithmic}[1]
\REQUIRE $A(n \times n, rowptr, colidx, val)$: csr matrix
\REQUIRE $ 0 < N_{b} <= n$: number of row-bands
\REQUIRE $ 1 < p <= N_{b}$: processor number ( or band number)
\STATE $numb\_cols\_per\_band \leftarrow n / N_{b}$
\STATE $N_p \leftarrow numb\_cols\_per\_band$
\IF{$p = N_{b} - 1$}
\STATE $N_p \leftarrow n - p \times numb\_cols\_per\_band$
\ENDIF
\STATE $col_{displ} \leftarrow i \times numb\_cols\_per\_band$
\STATE $nnz_{displ} \leftarrow A.rowptr[row_{displ}]$
\STATE \COMMENT{\bf Recover $nnz_{count}$ of the band}
\STATE \COMMENT{\bf Allocate band-row matrix $A_{band}(n \times N_p, nnz_{count})$}
\STATE $c \leftarrow 1$
\FOR { $i=1$ \TO $n$ }
  \STATE $A_{band}.rowptr[i] \leftarrow c$
  \FOR { $k=A_{band}.rowptr[i]$ \TO $A_{band}.rowptr[i+1]$}
    \IF{$nnz_{displ}<=A.colidx[k]<nnz_{displ}+N_p$}
      \STATE $A_{band}.colidx[c] \leftarrow A.colidx[k] - nnz_{displ}$
      \STATE $A_{band}.val[c] \leftarrow A.val[k]$
      \STATE $c \leftarrow c + 1$;
    \ENDIF
  \ENDFOR
\ENDFOR
\end{algorithmic}
\end{algorithm}

\subsection{Block-diagonal splitting}
\label{subsec:block-diagonal_splitting}

Table~\ref{tab:block-diagonal_partitioning} gives the schema of block-diagonal partitionning. In this section we briefly highlight the particularity of the diagonal block and his computation in parallel.
\begin{table}[!ht]
\SetDisplayIndex{1}
\SetSplitType{0}
\def\procSy{3}
\def\procEy{6}
\def\procSa{3}
\def\procEa{6}
\def\procSx{3}
\def\procEx{6}
\centering
\setlength{\tabcolsep}{0.0pt}
\begin{tabular}{cccc}
$\;\;\;\;$ & $\;\;\;\;\;\;\;=$ & $\;$ & $\;\;\;\;$ \\
\scalebox{0.45}{
\def\nnzcoef{1/1,2/1,3/1,4/1,5/1,6/1,7/1,8/1,9/1,10/1}
\SetDisplayInfoSize{0}
\SetDisplayInfoSizeProc{0}
\DrawBandDiagonalVector{(0,0)}{10}{1}{\nnzcoef}}
& &
\scalebox{0.45}{
\def\nnzcoef{
1/1,1/2,1/4,4/6,
2/1,2/2,2/3,4/5,6/8,
3/2,3/3,3/4,3/5,
4/1,4/3,4/4,9/10,8/7,
5/4,5/3,5/5,5/7,
6/4,6/6,6/7,
7/8,7/5,7/6,7/7,7/9,7/10,
8/6,10/9,8/8,8/10,
9/9,
9/7,10/7,10/8,10/10}
\SetDisplayInfoSize{0}
\SetDisplayInfoSizeProc{0}
\DrawBandDiagonalMatrix{(0,0)}{10}{10}{\nnzcoef}}
&
\scalebox{0.45}{
\def\nnzcoef{1/1,2/1,3/1,4/1,5/1,6/1,7/1,8/1,9/1,10/1}
\SetDisplayInfoSize{0}
\SetDisplayInfoSizeProc{0}
\DrawBandDiagonalVector{(0,0)}{10}{1}{\nnzcoef}
}\\
$\;\;\;Y$ & $\;\;\;\;\;\;\;$ & $A$ & $\;\;\;X$
\end{tabular}
\caption{Block-diagonal partitionning}
\label{tab:block-diagonal_partitioning}
\end{table}
When the product is performed by the matrix, the product of the diagonal block requires only local terms of the vector $x$. In contrast, off-diagonal coefficients require the corresponding terms of the vector $x$. The diagonal block are thick black lines in Figure~\ref{tab:block-diagonal_partitioning}.

The optimal splitting is the one that partitions the mesh into sub-structures of the same size, in order to balance the load with a smallest possible boundary to limit data transfers. Sub-structures should be as spherical as possible topologically, since it is the sphere which has the smaller outer surface.

\subsection{Sub-structuring splitting}
\label{subsec:sub-structuring_splitting}

In order to illustrate the sub-structuring method we consider a problem coming from the finite element discretization of an elliptic partial differential problem. To simplify the analysis, we consider the Laplace equation. However, the analysis can be carry out to any coercive elliptic problem. The model problem for the unknown $u$, in a bounded domain $\Omega$ with homogeneous Dirichlet boundary conditions on the boundary $\partial \Omega = \Gamma$ can be expressed as: for $f\in L^2(\Omega)$, find $u\in H^1(\Omega)$ such that
\begin{eqnarray}
\label{eq:model_problem}
-\nabla^2 u = f & in\; \Omega \\
u = 0 & on\; \Gamma
\end{eqnarray}
An equivalent variational formulation of this problem can be formulated as: for $f\in L^2(\Omega)$, find $u\in H_0^1(\Omega)$ such that
\begin{eqnarray}
\label{eq:formulation_variationnel}
\forall v \in H_0^1(\Omega), & \int_\Omega \nabla u \nabla v = \int_\Omega fv
\end{eqnarray}
This problem is well posed, i.e., has one and only one solution. After a Galerkin discretization with finite elements and a choice of nodal basis, the linear system is obtained
\begin{equation}
\label{eq:linear_systems}
Su = f
\end{equation}
where $f$ denotes the right hand side, $x$ the unknown and $S$ the stiffness matrix which is a sparse, symmetric and positive-definite matrix. Conjugate gradient detailed in section~\ref{sec:conjugate_gradient} is used to solve this linear system.

\begin{figure}[!ht]
\centering
\includegraphics[scale=0.35]{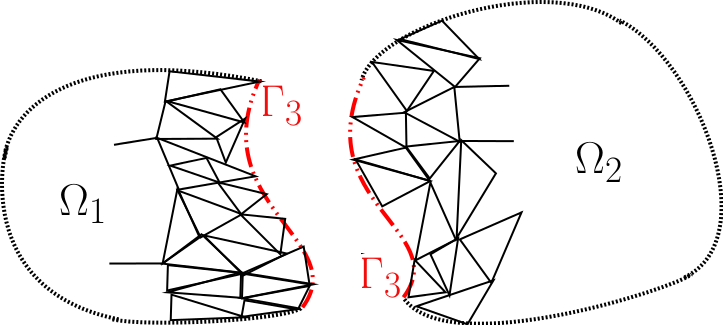}
\caption{Two sub-domains}
\label{fig:img:mesh_subdomains}
\end{figure}

In practice, mesh partitionning is a crucial step of finite element method. A finite element matrix is associated with a finite element mesh and the elements of the matrix are correlated with the interaction of the basis functions defined in the elements of the mesh. The total matrix is calculated as an assembly
of elementary matrices.
Let's consider a global domain $\Omega$ partitionned into two sub-domains without overlap $\Omega_1$ and $\Omega_2$, with a shared interface $\Gamma$ as drawn in Fig.~\ref{fig:img:mesh_subdomains}.
When a suitable numerotation of the degrees of freedom is harnessed, the stiffness matrix of the
initial considered model problem can be written as the following matrix:
\begin{eqnarray}
\label{eqn:2.4:eqn0}
  S =
  \left( \begin{array}{ccc}
  S_{11} & 0      & S_{13} \\
  0      & S_{22} & S_{23} \\
  S_{31} & S_{32} & S_{33}
  \end{array} \right)
\end{eqnarray}
It is formulated considering the case where the set of nodes numbered $1$ and $2$ are respectively associated to the sub-domains $\Omega_1$ and $\Omega_2$.
The last set of nodes numbered $3$ corresponds to the interface nodes of both sub-domains.
The corresponding linear system for previous matrix (\ref{eqn:2.4:eqn0}) can be reformulated as follows
\begin{eqnarray}
\label{eqn:2.4:eqn1}
  \left( \begin{array}{ccc}
  S_{11} & 0      & S_{13} \\
  0      & S_{22} & S_{23} \\
  S_{31} & S_{32} & S_{33}
  \end{array} \right)
  \left( \begin{array}{ccc}
  x_1 \\
  x_2 \\
  x_3
  \end{array} \right)
  =
  \left( \begin{array}{ccc}
  f_1 \\
  f_2 \\
  f_3
  \end{array} \right)
\end{eqnarray}
where $x=(x_1,x_2,x_3)^t$ is the unknown vector and $f=(f_1,f_2,f_3)^t$ is the right hand side.
The blocks $S_{13}$ and $S_{23}$ are respectively the transpose matrix of $S_{31}$ and $S_{32}$, and the blocks $S_{11}$ and $S_{22}$ are symmetric positive-definite whether $S$ was symmetric
positive-definite.
By assigning the different sub-domains at distinct processors, the local matrices can be formulated in parallel as follows:
\begin{eqnarray}
\label{eqn:2.4:mat_subdomains}
  S_{1} =
  \left( \begin{array}{cc}
  S_{11} & S_{13} \\
  S_{31} & S_{33}^{(1)}
  \end{array} \right)
  ~,~
  S_{2} =
  \left( \begin{array}{cc}
  S_{22} & S_{23} \\
  S_{32} & S_{33}^{(2)}
  \end{array} \right)
\end{eqnarray}
The blocks $S_{33}^{(1)}$ and $S_{33}^{(2)}$ denotes the interaction between the nodes on the interface $\Gamma$, respectively integrated in sub-domains $\Omega_1$ and on $\Omega_2$, i.e.,
\begin{equation}
\label{eq:interface_summation}
S_{33} = S_{33}^{(1)} + S_{33}^{(2)}
\end{equation}
In practice, the sub-domains $\Omega_1$ and $\Omega_2$ respectively knows the set of nodes ($1$, $3$) and ($2$,$3$)).
In order to avoid deadlock, the list of neighboring interface is re-ordering using the Welsh-Powell algorithm~\cite{welsh_upper_1967,brelaz_new_1979} for graph coloring.

\section{Conjugate Gradient (CG) for solving $Ax = b$}
\label{sec:conjugate_gradient}

Among all iterative methods, the Conjugate Gradient method (CG) is very efficient for solving the linear system $Ax=b$ when $A$ is a symmetric positive-definite matrix. CG algorithm has the advantage of being effective and easy to implement~\cite{{saad_iterative_2003}}.
The algorithm is based on the minimization of the distance of the iterative solution ($||Ax-b||_2$) to the exact solution on Krylov subspaces.

If we multiply $S$ by $v$ to obtain $Sv$, we can then multiply $S$ by the last obtained vector to get $S^2v$, and etc., so it is trivial to construct a Krylov sequences
$$v, Sv, S^2v, S^3v, S^4v, S^5v, ...$$
Krylov subspace methods are a large category of iterative methods used to solve linear systems from a wide range of applications. As introduced, in each iteration one or more sparse matrix-vector products is used, and Krylov subspace methods add vector(s) to a basis for one or more Krylov subspace(s)
\begin{equation}
\label{eq:krylov_subspace}
V^{(s)} = Span\{v_{1}, Sv_{1},\ldots,S^{s-1}v_{1}\}
\end{equation}
where $v_{1}$ is the initial vector solution, which is equal to the initial residual $r_1=f-Sx_1$.

\subsection{General conjugate gradient algorithm}
\label{subsec:preconditioned_CG_algorithm}

Knowing how data are distributed on each processors, we now present the conjugate gradient method for each type of distribution.
Algorithm~\ref{algo:cg-general_algorithm} gives the basic conjugate gradient algorithm with preconditionner. Diagonal preconditioner is considered in this paper. We denote $x_1$ the initial vector solution, $r_k$ the residual, $x_k$ the solution, $p_k$ the descent direction vector at the iteration $k$ and $(.,.)$ is the Euclidean scalar product. At each iteration, a sparse matrix-vector product, which is the most time consuming operation, has to be performed. The other operations are only scalar products and linear combinations of vectors.
\begin{algorithm}
\footnotesize
\caption{Conjugate Gradient method}
\label{algo:cg-general_algorithm}
\begin{algorithmic}[1]
\STATE Compute $r_1=f-Sx_1$, $z_1=M^{-1}r_1$, set $p_1=z_1$%
\FOR{$k=1,2,\ldots$} %
  \STATE $\alpha_k=(z_k,r_k)/(Sp_k,p_k)$%
  \STATE $x_{k+1}=x_k+\alpha_k p_k$%
  \STATE $r_{k+1}=r_k-\alpha_k Sp_k$%
  \STATE $z_{k+1}=M^{-1}r_{k+1}$%
  \STATE $\beta_k=(r_{k+1},z_{k+1})/(r_k,z_k)$
  \STATE $p_{k+1}=z_{k+1}+\beta_k p_k$%
\ENDFOR
\end{algorithmic}
\end{algorithm}
The CPU and GPU code are similar, as described in Algorithm~\ref{algo:cg-general_algorithm}, except that in GPU version the linear algebra operations are performed on GPU.
Many studies ~\cite{gaikwad_parallel_2010,GPU:BCK:2011,GPU:MMKPJBL:2013,cheikahamed:2012:inproceedings-1,cheikahamed:2013:inproceedings-4,ahamed2013stochastic,GPU:LS:2010,GPU:LSOLS:2010} in sparse linear systems demonstrated the effectiveness of GPU Computing compared to sequential and parallel CPU code for large size matrices. The performance are more better when memory~\cite{ref:cheikahamed:2012:inproceedings-1:DARG:2013,ref:cheikahamed:2012:inproceedings-1:DARG:2013:b,creel_high_2012} is well managed and threading distribution is well tuned~\cite{GPU:DZO:2011,cheikahamed:2012:inproceedings-2}. In this paper, we use Alinea, our research group library that offers linear algebra operations in both CPU and GPU CUDA 4.0~\cite{GPU:CUDA4.0:2011,GPU:BYFWA:2009}. The implementation of this library are analyzed in~\cite{cheikahamed:2012:inproceedings-1,cheikahamed:2012:inproceedings-2,cheikahamed:2013:inproceedings-3,cheikahamed:2013:inproceedings-4} and have shown its effectiveness and robustness compared to Cusp\cite{GPU:CUSP0.3.0:2012}, CUBLAS~\cite{GPU:CUBLAS4.0:2011}, CUSPARSE~\cite{GPU:CUSPARSE4.0:2011} for double precision arithmetics.

\subsection{Band-row algorithm}
\label{subsubsec:band-row_algorithm}

\subsubsection{Sparse matrix-vector product}
\label{subsubsec:band-row_spmv}

The processor that will perform the matrix-vector product for a band-row has only the corresponding terms of the vector $x$, colored area in Table~\ref{tab:band-row_splitting}. In order to carry out the sparse matrix-vector, this process needs all the terms of the vector $x$. The first step consists therefore to collect the terms that lacks, located out of the colored area in Figure 2. As it is the same for all processors, it will therefore be necessary to reconstruct the full vector $x$ on each processor. This operation corresponds to a classic collective exchange, where everyone is both transmitter and receiver.

\begin{figure}[!ht]
\centering
\includegraphics[scale=0.275]{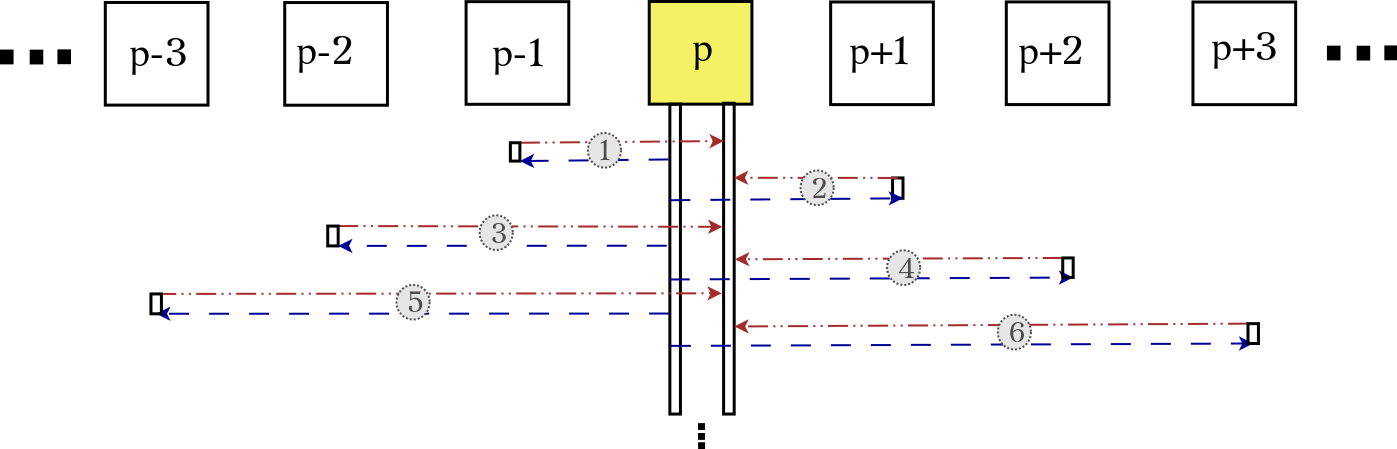}
\caption{Send/Recv ordering f processor number $p$}
\label{fig:img:send_recv_ordering}
\end{figure}

In this work, instead of using the collective operation, \emph{MPI\_Allgather}, including in message passing library (MPI), we use the equivalent \emph{Send/Recv}, with a \emph{left-right} ordering of sending and receiving. For the processor $p$, the \emph{left-right} ordering consists in respectively sending and receiving to and from $k=p-1$, $k=p+1$, $k=p-2$, $k=p+2$, $k=p-3$, $k=p+3$, ..., if $k > 0$.  This process is described in Fig.~\ref{fig:img:send_recv_ordering}.
The number of arithmetical operations requires to perform the local sparse matrix-vector multiplication, which is approximately $\frac{K \times n}{s}$, where $s$ is the number of processors, $n$ the dimension of the matrix, and $K$ the average number of non-zero coefficients per row.
On the other hand, the total number of terms of the vector $x$ to recover before performing the product is approximately $\frac{(s-1).n}{s}$, if the local matrix has non-zero values in almost all columns.
The amount of data is not small compared with the number of arithmetic operations. To optimize communications, it consists in finding a way to limit drastically the number of external values of vector $x$, located on the others processors, and necessary to compute the product by the matrix.

\subsubsection{Basic linear algebra operations}
\label{subsubsec:band-row_blalgebraop}

The computation of dot product is a relatively simple operation. Each processor performs local dot product, i.e., multiplies its elements and sums them, from their two local vectors. Finally, the local sums are added using \emph{MPI\_Allreduce} with \emph{MPI\_SUM} operation. Then, each processor has the global dot product. The operations such as addition of vectors, element wise product, etc. do not change compared to the sequential code. For GPU version, local operations have been performed on graphics card.

\subsection{Band-column algorithm}
\label{subsubsec:band-column_algorithm}

\subsubsection{Sparse matrix-vector product}
\label{subsubsec:band-column_spmv}

Unlike the band-row splitting sparse matrix-vector multiplication, the SpMV for band-column splitting avoid the exchange of the vector $x$. However, an $MPI\_Allreduce (MPI\_SUM)$ is required, in order to assemble the vector $y = Ax$. Note that basic operations are the same as for band-row splitting. For GPU version, local operations have been performed on graphics card. For the matrix-vector product, the results are first send to CPU before applying the same procedure as in CPU. 

\subsection{Sub-structuring method}
\label{subsubsec:ss_algorithm}

\subsubsection{Sparse matrix-vector product}
\label{subsubsec:sub-structuring_spmv}

Knowing how matrix is partitionned into sub-structures, Section~\ref{subsec:sub-structuring_splitting} now focus on the analysis of the sub-structuring method based on CG algorithm described in Algorithm~\ref{algo:cg-general_algorithm}.
This algorithm requires to perform a multiplication of the matrix $S$ by a descent direction vector $x=(x_1,x_2,x_3)^t$ at each iteration.
With the splitting into two sub-domains, the global matrix-vector multiplication can be written as follows:
\begin{eqnarray*}
\label{eqn:2.4:spmv}
  \left( \begin{array}{ccc}
  y_1 \\
  y_2 \\
  y_3
  \end{array} \right)
  & = &
  \left( \begin{array}{ccc}
  S_{11} & 0      & S_{13} \\
  0      & S_{22} & S_{23} \\
  S_{31} & S_{32} & S_{33}
  \end{array} \right)
  \left( \begin{array}{ccc}
  x_1 \\
  x_2 \\
  x_3
  \end{array} \right) \\
  & = &
  \left( \begin{array}{ccc}
  S_{11}x_1 + S_{13}x_3 \\
  S_{22}x_2 + S_{23}x_3 \\
  S_{31}x_1 + S_{32}x_2 + S_{33}x_3
  \end{array} \right)
\end{eqnarray*}
Considering the local matrices described equation~(\ref{eqn:2.4:mat_subdomains}), we can independently compute the both local matrix-vector products as follows 
\begin{eqnarray*}
  \left( \begin{array}{c}
  y_{1}\\
  y_{3}^{(1)}
  \end{array} \right)
  =
  \left( \begin{array}{cc}
  S_{11}x_1 + S_{13}x_3 \\
  S_{31}x_1 + S_{33}^{(1)}x_3
  \end{array} \right)
\end{eqnarray*}
\begin{eqnarray*}
  \left( \begin{array}{c}
  y_{2}\\
  y_{3}^{(2)}
  \end{array} \right)
  =
  \left( \begin{array}{cc}
  S_{22}x_2 + S_{23}x_3 \\
  S_{32}x_2 + S_{33}^{(2)}x_3
  \end{array} \right)
\end{eqnarray*}
Since $S_{33} = S_{33}^{(1)} + S_{33}^{(2)}$, and $y_{3} = y_{3}^{(1)} + y_{3}^{(2)}$. According to this last remark, SpMV can be calculated in two steps:
\begin{itemize}
  \item calculate the local matrix-vector multiplication in each sub-domain
  \item assemble on the interface, the local contributions
\end{itemize}
The first step involves only local data. The second requires the exchange of data between processes dealing with sub-domains with a common interface. In order to assemble interface values of neighboring sub-domains, each processor responsible to a sub-domain must know the description of its interfaces.

\subsubsection{Exchange at the interfaces}
\label{subsubsec:interf_exchange}

\begin{figure}[!ht]
\centering
\includegraphics[scale=0.45]{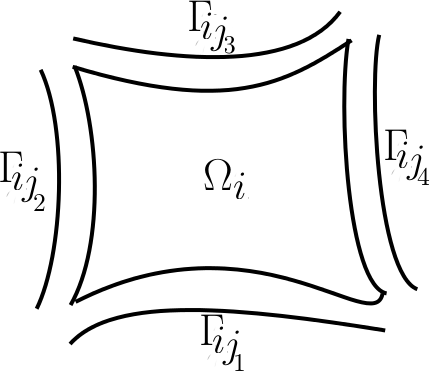}
\caption{Interface description}
\label{fig:img:domain_interfaces}
\end{figure}

When a sub-domain $\Omega_i$ has several neighboring sub-domains, we denote $\Gamma_{ij}$ the interface between $\Omega_i$ and $\Omega_j$ as described in Fig.~\ref{fig:img:domain_interfaces}.
An interface is identified by its neighboring sub-domains and the equations associated with its nodes. The interface is evaluated from its sub-domains using sparse matrix-vector product. This computation is in two steps for each neighboring sub-domain (Algorithm~\ref{alg:collect_send_neighbor}): \underline{\emph{collect}} the values of the local vector $y = Ax$ for all interfaces nodes, and then \underline{\emph{send}} this list to the $y$ vector of the interface equation.
\begin{algorithm}
\footnotesize
\caption{Construct inner buffer and send to neighboring}
\label{alg:collect_send_neighbor}
\begin{algorithmic}
\FOR{$s = 1$ \TO $number\_of\_neighboring$}
  \FOR{$i = 1$ \TO $n_s$}
    \STATE $buffer_s(i) = y(list_s(i))$
  \ENDFOR
  \STATE Send $buffer_s$ to neighbour(s)
\ENDFOR
\end{algorithmic}
\end{algorithm}
The next step consists in updating these changes to all neighboring sub-domains at interfaces equations. First, the contributions of the array containing the result of the matrix-vector product at the interface are received and then values on corresponding interface nodes are updated. This processus is described in Algorithm~\ref{alg:recv_update}. In GPU code, the construction of the inner buffer is carried out on CPU before sending it to the neighboring sub-domain.
\begin{algorithm}
\footnotesize
\caption{Receiving interface results and updating interface equations}
\label{alg:recv_update}
\begin{algorithmic}
\FOR{$s = 1$ \TO $number\_of\_neighboring$}
  \STATE Receive $buffer_s$ from neighbour(s)
  \FOR{$i = 1$ \TO $n_s$}
    \STATE $y(list_s(i)) = y(list_s(i)) + buffer_s(i)$
  \ENDFOR
\ENDFOR
\end{algorithmic}
\end{algorithm}
When an equation is shared by several interfaces,  the node value of the local vector $y = Ax$  in question is sent to all interfaces to which it belongs.
For any number of sub-domains, the mechanism of interface exchange and update is similar to those previously presented.
In GPU code, the procedure of exchange has been performed on CPU and then the assembled vector is copied back to GPU, before continuing the algorithm.

The use of the sub-structuring approach in the iterative GC algorithm is inherently parallel and making it excellent candidate for  implementation on parallel computers.
Indeed, we can distribute the sub-domains over all the available processors and thus compute locally the matrix-vector products, independently in parallel, and use distributed memory in order to limit the memory usage. As explained previously, after the computing of local matrix-vector multiplications, they required to be assembled along the interface. The key ingredient of the data is the local matrix $C$ that arise from the finite element discretization, with this approach, each node $i$ only needs to store $C_i$, the corresponding local matrix  to the sub-domain $\Omega_i$, which is only a fraction of the original matrix.
The dot product requires that each processor compute a weighted combination of the interface contributions in order to update their own data. After that, an $MPI\_Allreduce$ is required to compute the global inner product. Transfers between CPU and GPU at each iteration can decrease the performance of the exchanges algorithms. 

This new version of the CG algorithm based on sub-structuring method introduces only two new steps, which reside on data exchange. They consist firstly to share the contributions of local computed SpMV at the interface. Each machine requires to know the list of nodes along the interface and the number of neighboring sub-domains. Secondly, to assemble results over the cluster in order to piece together the local scalar product. This action, realized with MPI, is independant to the splitting.
Finally, another advantage of this algorithm is that it can easily be generelised for $n$ sub-domains.
This approach presents howewer two disavantages. The first drawback arises from an informatic point of view. The granularity, i.e., the number of operations to be performed by the processors compared to the amount of data received or send
by the processors may be weak. Indeed, here the granularity is proportional to the number of nodes in the sub-domains compared to the number of nodes on the interface.
The number of operations depend on the first parameter and the data transfer depends on the second parameter. If a lot of sub-domains are used, the interface size will not be small compared to the local sub-problem size. This means that the processors realize few computing operations (a local matrix-vector product) and a lot of communications. The second and more important drawback is an algorithmic one.
The classical parallel preconditioners per sub-domain are based on incomplete factorization of the local matrices. Such preconditioners are less and less efficient when the number of sub-domain increases.
The efficiency due to the parallelism is balanced by a slow convergence of the iterative CG method.

\noindent {\sc Remark:} A way to avoid this problem, is to use some preconditioners based on domain decomposition methods~\cite{NME:NME1620320604,quarteroni_domain_1999,toselli_domain_2004,SBG1996,magoules:journal-auth:16,magoules:journal-auth:21}.
Optimized transmissions conditions between the sub-domains are usually defined like in~\cite{magoulesf_contrib_3:Chevalier:1998:SMO,magoulesf_contrib_3:Gander:2000:OSM,magoules:journal-auth:9,magoules:journal-auth:10,magoules:journal-auth:13,magoules:journal-auth:14,magoules:journal-auth:18,magoules:journal-auth:23} for the Schwarz method, or in~\cite{magoules:journal-auth:24} for the FETI method.
Optimization of the transmissions conditions based on small patch also exists as first introduced in~\cite{magoules:proceedings-auth:6,magoules:journal-auth:8,magoules:journal-auth:12,magoules:journal-auth:29} for acoustics and in~\cite{magoules:journal-auth:17,magoules:journal-auth:20} for linear elasticity.
The problem on the interface, is then solved on the CPU with an iterative algorithm.
At each iteration, each sub-problem defined in each sub-domain is solved on the GPU with the CG method.
This hybrid CPU-GPU domain decomposition method implementation was first proposed in~\cite{Papadrakakis:2011} for the FETI method and in~\cite{cheikahamed:2013:inproceedings-4}, for the optimized Schwarz method, but is outside the scope of this paper.
In the following we present the results obtained with the CG algorithm issued from a sub-structuring approach as described in the previous paragraphs.

\section{Experimental Evaluation}
\label{sec:exp_evaluation}

This section reports and presents the evaluation of the a set of numerical experiments.

\subsection{Hardware platform}
\label{subsec:experiments_platform}

For our performance analysis we considered a machine based on an Intel Core i7 920 2.67Ghz, which has 4 physical cores and 4 logical cores, 12GB RAM, and two different system accelerated generations of nVidia graphics card: GTX275 with 895MB memory, which is double precision compatible.

\subsection{Matrices tested}
\label{subsec:matrices_tested}

In order to evaluate our analysis on large scale engineering problems, we use a set of matrices from the University of Florida repository~\cite{davis_university_2011}.
The properties of these matrices are reported in Table~\ref{tab:sketches_matrices} where $h$, $nz$, $density$, $bandwidth$, $max row$, $nz/h$ and $nz/h stddev$ present the size of the matrix, the number of non-zero coefficients, the density corresponding to the number of non-zero values divided by the total number of matrix elements, the upper bandwidth, the maximum row density, the mean row density and the standard deviation of $nz/h$. A pattern of the structure matrix is collected in the first column, and the histograrow density is described in the second column.
\begin{table}
\centering
\scalebox{0.85}{
\begin{tabular}{m{6.4cm}r}
\hlinewd{1.0pt}
\infosparse{img/png/sparse/qa8fm}{qa8fm}{66127}{1660579}{0.038}{1048}{27}{25.112}{4.183}\\
\multicolumn{2}{c}{\emph{3D acoustic FE mass matrix. A. Cunningham, Vibro-Acoustic Sciences Inc.}}\\
\hlinewd{1.0pt}
\infosparse{img/png/sparse/2cubes_sphere}{2c\_sphere}{101492}{1647264}{0.016}{100407}{31}{16.230}{2.654}\\
\multicolumn{2}{c}{\emph{FEM, electromagnetics, 2cubes in a sphere. Evan Um, Geophysics, Stanford.}}\\
\hlinewd{1.0pt}
\infosparse{img/png/sparse/thermal2}{thermal2}{1228045}{8580313}{0.001}{1226000}{11}{6.987}{0.811}\\
\multicolumn{2}{c}{\emph{Unstructured FEM, steady state thermal problem. Dani Schmid, Univ. Oslo.}}\\
\hlinewd{1.0pt}
\infosparse{img/png/sparse/thermomech_TK}{thermo\_TK}{102158}{711558}{0.007}{102138}{10}{6.965}{0.715}\\
\multicolumn{2}{c}{\emph{FEM problem, temperature and deformation of a steel cylinder.}}\\
\hlinewd{1.0pt}
\infosparse{img/png/sparse/cfd2}{cfd2}{123440}{3085406}{0.020}{4332}{30}{24.995}{3.888}\\
\multicolumn{2}{c}{\emph{CFD, symmetric pressure matrix, from Ed Rothberg, Silicon Graphics, Inc.}}\\
\hlinewd{1.0pt}
\infosparse{img/png/sparse/finan512}{finan512}{74752}{596992}{0.011}{74724}{55}{7.986}{6.278}\\
\multicolumn{2}{c}{\emph{Economic problem}}\\
\hlinewd{1.0pt}
\infosparse{img/png/sparse/Dubcova2}{Dubcova2}{65025}{1030225}{0.024}{64820}{25}{15.844}{5.762}\\
\multicolumn{2}{c}{\emph{Univ. Texas at El Paso, from a PDE solver.}}\\
\hlinewd{1.0pt}
\infosparse{img/png/sparse/af_shell8}{af\_shell8}{504855}{17579155}{0.007}{4909}{40}{34.820}{1.285}\\
\multicolumn{2}{c}{\emph{Olaf Schenk, Univ. Basel: AutoForm Eng. GmbH, Zurich. sheet metal forming.}}\\
\hlinewd{1.0pt}
\end{tabular}
}
\caption{Sketches of finite element matrices}
\label{tab:sketches_matrices}
\end{table}

\subsection{Numerical results}
\label{subsec:numerical_results}

The time in seconds (s) of partitionning in row-bands, column-bands and sub-structures are reported in Table~\ref{tab:time_-partitionning}. The first column gives the name of the matrix. In the second column are collected sub-structuring splitting using the metis software at graph coloring step. The band column splitting are given in third column. The last column gives the band-row splitting.
\begin{table}
\centering
\begin{tabular}{lccc}
\hlinewd{1.0pt}
\textbf{Matrix} & \textbf{metis} & \textbf{band-c} & \textbf{band-r}\\
\hlinewd{1.0pt} %
\multicolumn{4}{c}{\bf \it 2-partitionning}\\
\texttt{qa8fm} & 3.16 & 1.15 &  1.20 \\
\texttt{2c\_sphere} & 5.40 & 2.49 & 2.49\\
\texttt{thermo\_TK} & 1.58 & 1.03 & 1.0\\
\texttt{cfd2} & 6.74 & 3.72 & 2.06\\
\texttt{thermal2} & 12.53 & 6.98 & 6.53\\
\texttt{af\_shell8} & 68.05 & 26.93 & 27.65\\
\texttt{finan512} & 0.83 & 0.44 & 0.42\\
\hlinewd{1.0pt}
\multicolumn{4}{c}{\bf \it 4-partitionning}\\
\texttt{qa8fm} & 4.60 & 1.19 &  1.11 \\
\texttt{2c\_sphere} & 5.14 & 2.52 & 2.49\\
\texttt{thermo\_TK} & 2.83 &1.03  & 1.0\\
\texttt{cfd2} & 8.28 & 2.52 & 2.04\\
\texttt{thermal2} & 15.72 & 7.37 & 7.63\\
\texttt{af\_shell8} & 86.18 & 28.32 & 51.36\\
\texttt{finan512} & 1.01 & 0.45 & 0.42\\
\hlinewd{1.0pt}
\multicolumn{4}{c}{\bf \it 8-partitionning}\\
\texttt{qa8fm} & 6.64 & 1.28 &  1.51 \\
\texttt{2c\_sphere} & 6.08 & 2.58 & 2.49\\
\texttt{thermo\_TK} & 1.88 & 1.07 & 1.0\\
\texttt{cfd2} & 10.84 & 2.25 & 2.61\\
\texttt{thermal2} & 17.33 & 8.01 & 7.75\\
\texttt{af\_shell8} & 104.41 & 28.72 & 27.38\\
\texttt{finan512} & 1.21 & 0.48 & 0.42\\
\hlinewd{1.0pt}
\end{tabular}
\caption{Execution time of partitionning (s)}
\label{tab:time_-partitionning}
\end{table}
Table~\ref{tab:parallel_cg_substructuring_result_dp_csr} reports respectively the running times of sub-structuring CG in seconds (s) (CSR) on CPU and GPU. The expermiments have been performed on a single node.
\begin{table}
\centering
\setlength{\tabcolsep}{1.5pt}
\begin{tabular}{llllllllll}
\hlinewd{1.0pt}
{Matrix} & \#iter. & 1CPU & 2CPU & {2CPUs} & {2GPUs} & {4CPUs} & {4GPUs} & {8CPUs} & {8GPUs} \\
\hlinewd{1.0pt} %
\texttt{2c\_sphere} & 24 & 0.386 & 0.026 & 0.209 & 0.047 & 0.124 & 0.065 & 0.125 & 0.113\\
\texttt{af\_shell8} & 2815 & 374.356 & 14.549 & 198.998 & 23.422 & 110.662 & 18.814 & 107.668 & 21.385\\
\texttt{cfd2} & 2818 & 71.078 & 3.73 & 38.114 & 5.904 & 21.657 & 8.28 & 22.111 & 12.186\\
\texttt{Dubcova2} & 168 & 1.776 & 0.128 & 0.926 & 0.364 & 0.54 & 0.405 & 0.85 & 0.64\\
\texttt{finan512} & 15 & 0.117 & 0.017 & 0.121 & 0.063 & 0.067 & 0.071 & 0.145 & 0.12\\
\texttt{qa8fm} & 29 & 0.418 & 0.023 & 0.235 & 0.099 & 0.198 & 0.137 & 0.168 & 0.115\\
\texttt{thermo\_TK} & 13226 & 141.359 & 13.214 & 74.423 & 22.877 & 41.979 & 32.587 & 40.888 & 51.528\\
\hlinewd{1.0pt}
\end{tabular}
\caption{Execution time for parallel sub-structuring CG (s) for CSR format}
\label{tab:parallel_cg_substructuring_result_dp_csr}
\end{table}
Numerical results presented in Table~\ref{tab:parallel_cg_substructuring_result_dp_csr} clearly highlight the efficiency of GPU device compared to CPU parallel computation for solving linear systems with sub-structuring methods. The speed up to 19.9x for 4GPUs and 3.3x for 4CPUs (af\_shell8). The ratio between sequential CPU and 1GPU reaches 25.7x for the same matrix.

\section{Conclusion}
\label{sec:conclusion}

This paper gives an analysis of a parallel sub-structuring method based on conjugate gradient method for solving large and sparse linear systems on a cluster of GPU Computing. We have harnessed the efficiency of parallel algorithms, coupled with the high power computing of GPU. We have evaluated a parallelized conjugate gradient algorithm using sub-structuring method, which has a natural approach of parallelization.
The experiments have been performed on large sparse matrices arising from large scale engineering problems. The results clearly show the interest of sub-structuring accelerated with GPU Computing to solve linear systems for symmetric positive-definite matrices. The relative gains of the GPU cluster reaches 19x for 4GPUs compared to a 4 CPUs, and up to 27x compared to sequential CPU. 




\section*{Acknowledgment}
The authors acknowledge the CUDA Research Center at Ecole Centrale Paris (France) for its support and for providing the computing facilities.

\bibliography{bib/hpcc2014_paper1_ac_fm,bib/MAGOULES-JOURNAL1,bib/MAGOULES-PROCEEDINGS1}
\bibliographystyle{abbrv}

\end{document}